\documentclass{article}

\usepackage[utf8]{inputenc}
\usepackage[T1]{fontenc}
\usepackage{lmodern}
\usepackage{babel}
\usepackage{lipsum}
\usepackage{graphicx}
\usepackage{scalerel,stackengine,amsmath}
\usepackage{amsfonts}
\usepackage{subcaption}
\usepackage{tikz}
\usetikzlibrary{positioning, calc}
\usepackage[ruled,linesnumbered]{algorithm2e}
\usepackage{accents}
\usepackage{xcolor}
\usepackage{parskip}
\usepackage{amssymb}
\usepackage{hyperref}
\usepackage{multirow}
\usepackage{alphabeta} % Allows \alpha, etc. in text mode.
\usepackage{thm-restate}
\usepackage{listings}
\usepackage{tcolorbox}
\usepackage{tablefootnote}

\newif\ifneurips
\neuripsfalse

\newif\ifarxiv
\arxivtrue

\newif\ifopus
\opusfalse

\newif\ifnonjournal

\ifarxiv
    \nonjournaltrue
\else
    \ifopus
        \nonjournaltrue
    \else
        \nonjournalfalse
    \fi
\fi

\ifneurips
     \PassOptionsToPackage{numbers, compress}{natbib}
\usepackage{neurips_2024}

\usepackage[utf8]{inputenc} % allow utf-8 input
\usepackage[T1]{fontenc}    % use 8-bit T1 fonts
\usepackage{hyperref}       % hyperlinks
\usepackage{url}            % simple URL typesetting
\usepackage{booktabs}       % professional-quality tables
\usepackage{amsfonts}       % blackboard math symbols
\usepackage{nicefrac}       % compact symbols for 1/2, etc.
\usepackage{microtype}      % microtypography
\usepackage{xcolor}         % colors

\else

\ifarxiv
\usepackage{arxiv}
\usepackage{charter}
\usepackage{fullpage}
\fi

\ifopus
\usepackage{zibtitlepage}
\usepackage{arxiv}
\usepackage{charter}
\usepackage{fullpage}
\fi

\hypersetup{
    colorlinks=true,
    linkcolor=blue,
    filecolor=magenta,      
    urlcolor=cyan,
    citecolor=blue
    }

\fi

\newcommand{\scipml}{PySCIPOpt-ML\xspace}
\newcommand{\reals}{\ensuremath{\mathbb{R}}\xspace}

\ifneurips
\author{%
  Mark Turner \\
  Department of Mathematical Optimization\\
  Zuse Institute Berlin\\
  Takustr. 7, 14195 Berlin\\
  \texttt{turner@zib.de} \\
  Antonia Chmiela \\
  Department of Mathematical Optimization \\
  Zuse Institute Berlin\\
  Takustr. 7, 14195 Berlin\\
  \texttt{chmiela@zib.de} \\
  Thorsten Koch \\
  Institute of Mathematics \\
  Technische Universit{\"a}t Berlin \\
  Straße des 17. Juni 135, 10623 Berlin \\
  \texttt{koch@zib.de} \\
  Michael Winkler \\
  Gurobi GmbH \\
  Ulmenstr. 37-39, 60325 Frankfurt am Main \\
  \texttt{winkler@gurobi.com}
}
\else
\ifnonjournal
\author{ \href{https://orcid.org/0000-0001-7270-1496}{\includegraphics[scale=0.06]{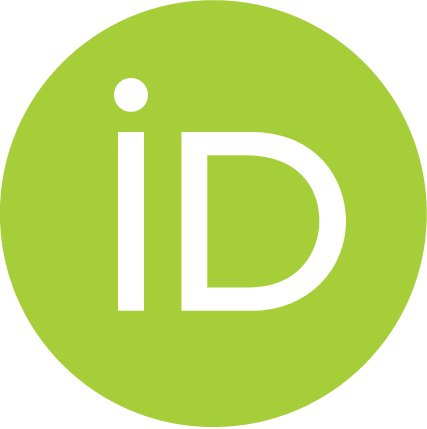}}\hspace{1mm}Mark Turner\thanks{Zuse Institute Berlin, Department of Mathematical Optimization, Takustr. 7, 14195 Berlin} \\
	\texttt{turner@zib.de} \\
	\And
	\href{https://orcid.org/0000-0002-4809-2958}{\includegraphics[scale=0.06]{orcid_id_icon.eps}}\hspace{1mm}Antonia Chmiela\footnotemark[1]\\
	\texttt{chmiela@zib.de} \\
	\And
	\href{https://orcid.org/0000-0002-1967-0077}{\includegraphics[scale=0.06]{orcid_id_icon.eps}}\hspace{1mm}Thorsten Koch\thanks{Chair of Software and Algorithms for Discrete Optimization, Institute of Mathematics, Technische Universit{\"a}t Berlin, Straße des 17. Juni 135, 10623 Berlin, Germany}\hspace{2mm}\footnotemark[1]
	\\
	\texttt{koch@zib.de}
    \And
	Michael Winkler\thanks{Gurobi GmbH, Ulmenstr. 37-39, 60325 Frankfurt am Main, Germany}\hspace{2mm} \\
	\texttt{winkler@gurobi.com}
}
\else
\author{Mark Turner\inst{1}\orcidID{0000-0001-7270-1496} \and
Antonia Chmiela\inst{1}\orcidID{0000-0002-4809-2958}\and
Thorsten Koch\inst{2,1}\orcidID{0000-0002-1967-0077}\and
Michael Winkler\inst{3}
}
\authorrunning{M. Turner et al.}
\titlerunning{PySCIPOpt-ML}
\institute{Zuse Institute Berlin, Germany\and
Institute of Mathematics, Technische Universit{\"a}t Berlin, Germany\and
Gurobi GmbH, Ulmenstr. 37-39, 60325 Frankfurt am Main, Germany\\
\email{\{turner, chmiela, koch\}@zib.de, winkler@gurobi.com}
}
\fi
\fi

\ifopus
\ZTPAuthor{\ZTPHasOrcid{Mark Turner}{0000-0001-7270-1496},
\ZTPHasOrcid{Antonia Chmiela}{0000-0002-4809-2958},
\ZTPHasOrcid{Thorsten Koch }{0000-0002-1967-0077}, \\
Michael Winkler}
\ZTPTitle{PySCIPOpt-ML: Embedding Trained Machine Learning Models into Mixed-Integer Programs}
\ZTPNumber{23-28}
\ZTPMonth{December}
\ZTPYear{2023}
\fi

\title{PySCIPOpt-ML: Embedding Trained Machine Learning Models into Mixed-Integer Programs}

\begin{document}

\ifopus
\zibtitlepage
\fi

\date{}
\maketitle

\begin{abstract}
    A standard tool for modelling real-world optimisation problems is mixed-integer programming (MIP). However, for many of these problems, information about the relationships between variables is either incomplete or highly complex, making it difficult or even impossible to model the problem directly. To overcome these hurdles, machine learning (ML) predictors are often used to represent these relationships and are then embedded in the MIP as surrogate models. Due to the large amount of available ML frameworks and the complexity of many ML predictors, formulating such predictors into MIPs is a highly non-trivial task. In this paper, we introduce PySCIPOpt-ML, an open-source tool for the automatic formulation and embedding of trained ML predictors into MIPs. By directly interfacing with a broad range of commonly used ML frameworks and an open-source MIP solver, PySCIPOpt-ML provides a way to easily integrate ML constraints into optimisation problems. Alongside PySCIPOpt-ML, we introduce, SurrogateLIB, a library of MIP instances with embedded ML constraints, and present computational results over SurrogateLIB, providing intuition on the scale of ML predictors that can be practically embedded.
    The project is available at \url{https://github.com/Opt-Mucca/PySCIPOpt-ML}.
\end{abstract}

\section{Introduction and Related Work}

Many problems coming from real-world applications can be modelled using \emph{Mixed-Integer Programming (MIP)}, which can be expressed as follows:
\begin{equation} \label{mip} \tag{P}
\begin{aligned}
& \underset{x \in \mathbb{R}^n}{\text{min}}
& & f_0(x) \\
& \text{s.t.}
& & f_i(x) \leq 0, \ \forall i \in [m], \\
&&& x_i \in \mathbb{Z}, \ \forall i \in \mathcal{I}.
\end{aligned}
\end{equation}
Each function $f_i: \mathbb{R}^n \rightarrow \mathbb{R}$ for $i \in [m]$ is continuous, and $\mathcal{I} \subseteq [n]$ denotes the index set of the integer-constrained variables. Due to a MIP's ability to represent complex systems and decision-making processes across various industries, MIP has become a standard tool to solve optimisation problems arising in real-world applications \cite{achterberg2007constraint,petropoulos2023operational}.

A common challenge of finding a suitable MIP formulation \eqref{mip} of a system is the presence of unknown or highly complex relationships. This poses a problem since traditional MIP approaches rely on precise definitions of constraints and the objective function. To address this challenge, \emph{Machine Learning (ML)} predictors are often used to approximate these relationships and are then embedded in \eqref{mip} as surrogate models. To embed a trained ML predictor $g: \mathbb{R}^n \rightarrow \mathbb{R}$, it must first be formulated as a series of MIP constraints before it can be added to \eqref{mip}. However, with the growing popularity and accessibility of different ML frameworks \cite{tensorflow2015-whitepaper,scikit-learn,xgboost,ke2017lightgbm,paszke2019}, automatically formulating trained ML predictors coming from these different architectures into MIPs has become a highly non-trivial, yet necessary task.

\begin{figure}[t]
\centering
\includegraphics[width=\textwidth]{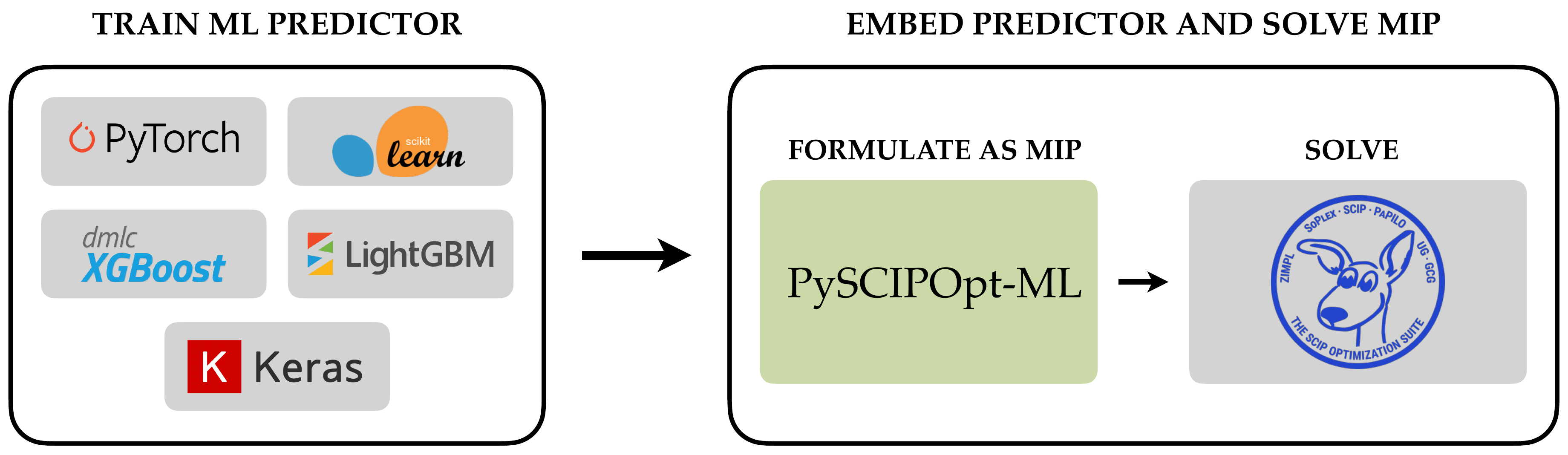}
\caption{How trained ML predictors can be automatically formulated as MIPs and embedded in an optimisation problem using our Python package \scipml.}
\label{figure}
\end{figure}

\paragraph{Contribution.} In this paper, we introduce \scipml, an open-source tool for the automatic formulation and embedding of trained ML predictors into MIPs. By directly interfacing with a broad range of commonly used ML frameworks, \scipml{} allows for the easy integration of ML constraints into MIPs, which can then be solved by the state-of-the-art open-source solver SCIP \cite{scip}. The package supports various models from Scikit-Learn \cite{scikit-learn}, XGBoost \cite{xgboost}, LightGBM \cite{ke2017lightgbm}, Keras \cite{keras}, and PyTorch \cite{paszke2019}, making it a general tool to easily optimise MIPs with embedded ML constraints (see Figure~\ref{figure}). To represent predictors as MIPs, we use formulations based on SOS1 constraints \cite{fischer2018branch}, which admit valid formulations in the absence of user specified bounds. Alongside \scipml, we introduce SurrogateLIB, an expandable library of MIP instances with embedded ML constraints based on real-world data. Computational results over this library of instances are presented, aiming to provide intuition on the scale of ML predictors that can be practically embedded.

To summarise, our contributions are the following:

\begin{enumerate}
    \item We introduce the \textbf{Python package \scipml{}
    \footnote{\url{https://github.com/Opt-Mucca/PySCIPOpt-ML}}},
    %\footnote{The code is provided in the supplementary materials.}},
    which \textbf{directly interfaces various ML frameworks with the open-source solver SCIP} (see Section~\ref{section:package}).
    \item We utilise \textbf{MIP formulations prioritising practicality}, since they do not rely on user defined bounds or expensive bound propagation techniques (see Section~\ref{section:formulation})
    \item We introduce a \textbf{new library of MIP instances with embedded ML constraints called SurrogateLIB}
    \footnote{\url{https://zenodo.org/records/10357875}},
    %\footnote{A compressed version of the library is provided in the supplementary materials. To respect the size limit of the submission, we had to remove instances.},
    as well as the corresponding instance generators (see Section~\ref{section:library}).
    \item We present \textbf{computational experiments showcasing the scale of predictors} that can be practically embedded (see Section~\ref{section:experiments}).
\end{enumerate}

\paragraph{Related Work.} In recent years, more and more industries have leveraged ML frameworks to approximate unknown or complex relationships in optimisation problems. For instance, ML predictors have been used to optimise cancer treatment plans \cite{bertsimas2016}, minimise energy consumption \cite{misaghian2023}, and in the design of energy systems \cite{JALVING2023121767}. Embedded ML predictors have also been used in the optimisation of gas production and water management systems \cite{lopez2023}, as well as in the optimisation of retail prices \cite{Ferreira2015} and locations \cite{huang2018}, to name only a few examples. For surveys of embedded ML predictors through surrogate optimisation, see \cite{lombardi2018boosting,bhosekar2018advances,mcbride2019overview}.

The increasing demand for embedding ML predictors into MIPs has motivated the development of various tools for their automatic integration. The commercial MIP solver Gurobi \cite{gurobi} released the Python package Gurobi Machine Learning \cite{gurobiml}, which integrates trained regression models into optimisation problems that can then be solved with Gurobi. Other packages that rely on Gurobi are OptiCL \cite{opticl}, which supports many Scikit-Learn predictors \cite{scikit-learn}, and Janos \cite{bergman2022janos}, which supports linear and logistic regression, as well as neural networks with rectified linear unit (ReLU) activation functions. EML \cite{eml} in contrast uses CPLEX \cite{cplex}, and enables the embedding of decision trees and neural networks through the Keras API \cite{keras}. ENTMOOT \cite{thebelt2021entmoot} focuses on representing gradient-boosted tree models from LightGBM \cite{ke2017lightgbm}. On the other hand, reluMIP \cite{reluMIP.2021} allows the user to integrate ReLU neural networks trained with TensorFlow \cite{tensorflow2015-whitepaper} within optimisation problems modelled with either Pyomo \cite{bynum2021pyomo} or Gurobi. Finally, OMLT \cite{ceccon2022omlt,ammari2023,zhang2024augmenting} focuses on embedding a wide variety of neural networks and gradient-boosted trees via the general MIP modelling framework Pyomo \cite{bynum2021pyomo} and the general purpose ML interface ONNX\footnote{\url{https://github.com/onnx/onnx}}.

While many of these packages focus on supporting specific types of predictors, \scipml{} allows the user to embed a wide variety of predictors. Specifically, for both regression and classification tasks: neural networks with various activation functions, linear and logistic regression models, decision trees, gradient boosted trees, random forests, support vector machines, and centroid clustering. In addition, unlike the other packages, \scipml interfaces directly with an open-source MIP solver. This allows users to dynamically and transparently improve the solving process, e.g., by dynamic bound tightening \cite{hojny2024verifying}.

\section{About the Package}
\label{section:package}

\scipml{} directly interfaces with the popular ML frameworks Scikit-Learn \cite{scikit-learn}, XGBoost \cite{xgboost}, LightGBM \cite{ke2017lightgbm}, Keras \cite{keras}, and PyTorch \cite{paszke2019}. Each framework allows the user to train ML predictors of different types, all of which can then be embedded into a MIP by our package. A detailed overview of what \scipml{} supports is given in Table~\ref{tab:ml_framework_overview}. 

\begin{table}[t]
\centering
%\scriptsize
%\setlength{\tabcolsep}{8pt} % Adjust the space between columns
%\renewcommand{\arraystretch}{1.3} % Adjust the space between rows
\resizebox{\columnwidth}{!}{%
\begin{tabular}{l|ccccc}
\hline
\hline \\
ML Predictor & PyTorch & Tensorflow / Keras & Scikit-Learn & LightGBM & XGBoost \\[0.5cm]
\hline
\hline \\
Decision Trees          & - & - &
\begin{tabular}{@{}c@{}} \textbf{SCIP}/EML/GRB/\\OCL/OMLT\end{tabular} &
\begin{tabular}{@{}c@{}} \textbf{SCIP}/EML/GRB/\\OMLT\end{tabular} &
\textbf{SCIP}/GRB/OMLT \\
&&&&&\\

Gradient Boosted Trees  & - & - &
\begin{tabular}{@{}c@{}} \textbf{SCIP}/GRB/OCL/\\OMLT\end{tabular} &
\begin{tabular}{@{}c@{}} \textbf{SCIP}/ENT/GRB/\\OMLT\end{tabular} &
\textbf{SCIP}/GRB/OMLT \\
&&&&&\\

Random Forests          & - & - &
\begin{tabular}{@{}c@{}} \textbf{SCIP}/GRB/OCL/\\OMLT\end{tabular} &
\begin{tabular}{@{}c@{}} \textbf{SCIP}/ENT/GRB/\\OMLT\end{tabular} &
\textbf{SCIP}/GRB/OMLT \\
&&&&&\\

\hline \\
Neural Networks (ReLU)        &
\textbf{SCIP}/GRB/OMLT &
\begin{tabular}{@{}c@{}} \textbf{SCIP}/EML/GRB/\\OMLT/RMIP\end{tabular} &
\begin{tabular}{@{}c@{}} \textbf{SCIP}/GRB/JAN/\\OCL/OMLT\end{tabular} &
- & - \\
&&&&&\\

Neural Networks (other)        & \textbf{SCIP}/OMLT & \textbf{SCIP}/EML/OMLT & \textbf{SCIP}/OMLT & - & - \\
&&&&&\\

\hline \\
Linear Regression       & - & - &
\begin{tabular}{@{}c@{}} \textbf{SCIP}/GRB/JAN/\\OCL\end{tabular} &
- & - \\
&&&&&\\

Logistic Regression     & - & - &
\begin{tabular}{@{}c@{}} \textbf{SCIP}/GRB/JAN/\\OCL\end{tabular} &
- & - \\
&&&&&\\

\hline \\
Support Vector Machines   & - & - &\textbf{SCIP}/OCL & - & - \\
&&&&&\\

Centroid Clustering   & - & - & \textbf{SCIP} & - & - \\
&&&&&\\

\hline
\hline
\end{tabular}
}
\caption{Overview of ML predictors supported by each automatic embedding framework. The following abbreviations are used: SCIP is \scipml, EML is Emperical Machine Learning, ENT is ENTMOOT, GRB is Gurobi Machine Learning, JAN is Janos, OCL is OptiCL, OMLT is OMLT, and RMIP is reluMIP. Predictor-framework pairs that do not exist or are not supported by any tool are indicated by "-".\\
\footnotesize{*For all supported frameworks and models \scipml offers multi-regression and multi-classification.}
}
\label{tab:ml_framework_overview}
\end{table}

\begin{figure}[t]
\centering
\includegraphics[width=\textwidth]{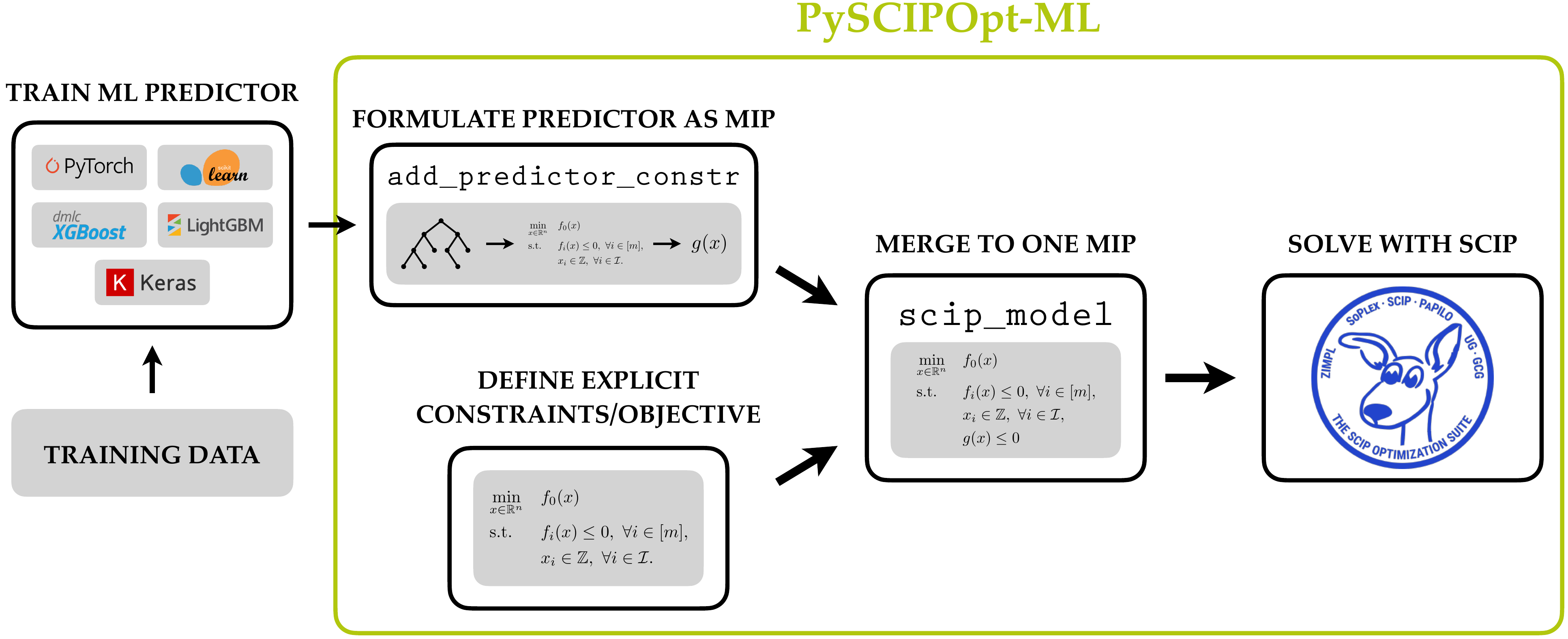}
\caption{The workflow when embedding ML predictors into a MIP with \scipml{}.}
\label{figure2}
\end{figure}

An advantage of \scipml{} is that it seamlessly joins all parts of the modelling, training, and optimisation processes. At no point is there a need to export the ML predictor or MIP model into some intermediary language. A central function (\texttt{add\_predictor\_constr}) acts as the main interface between the ML frameworks training the predictor and the solver SCIP solving the problem. To create and solve a MIP problem with an embedded ML predictor, the user must only do the following:
\begin{enumerate}
    \item Train the ML predictor, \texttt{predictor}, using one of the supported ML frameworks. 
    \item Define the optimisation problem, \texttt{scip\_model}. This involves adding the constraints and objective not directly related to the ML predictor. 
    \item Create variables that represent \texttt{(input\_data,output\_data)}. The output variables are optional as they can be inferred by the ML predictor. Embed the predictor into the MIP by simply calling:
    \[
    \small{\texttt{add\_predictor\_constr(scip\_model, predictor, input\_data, output\_data)}}.
    \]
    \item Solve the updated \texttt{scip\_model} with SCIP.
\end{enumerate}
Figure~\ref{figure2} visualises this workflow. 

In summary, \scipml{} helps to easily create, embed, and solve optimisation problems with ML surrogate models and spares the user from navigating multiple platforms at the same time.

\section{MIP Formulations of ML Predictors}
\label{section:formulation}

To embed an ML predictor in a MIP, we first need to represent the predictor itself as a MIP. The formulation of the predictor has a significant effect on the ability to efficiently optimise the larger optimisation problem. There has been a substantial amount of work on determining the best performing formulations for popular ML predictors. In particular, neural networks with ReLU activation functions \cite{grimstad2019relu,schweidtmann2019deterministic,anderson2020strong,kronqvist2021between} and ensemble based methods that use decision trees \cite{mistry2021mixed,thebelt2021entmoot,ammari2023,gurobiml}, have received a lot of attention in the past. Often, the most natural formulation of these predictors involves quadratic constraints. Due to performance issues, however, formulations that utilise big-M, and indicator constraints are in general preferred. In the remainder of this section, we will give an overview of how \scipml{} models different ML predictors as MIPs.

\subsection{Neural Networks}\label{subsec:nn-formulations}

The current standard MIP formulation for feed forward neural networks with ReLU activation functions uses big-M constraints \cite{grimstad2019relu,anderson2020strong,ceccon2022omlt,huchette2023deep, badilla2023computational}. Consider a single neuron in a fully connected hidden-layer $k$, where $k \in \mathbb{N}$. Let $x \in \mathbb{R}^{n}$ be the output of the previous layer $k-1$, $y \in \mathbb{R}_{\geq 0}$ be the output of the given neuron, and $w \in \mathbb{R}^{n}$ and $b \in \mathbb{R}$ be the weights and bias connecting the previous layer to the neuron. The standard big-M formulation requires bounds $l,u \in \mathbb{R}$ on the neuron's output, and a binary variable for the neuron, namely $z \in \{0, 1\}$. The formulation is defined as:
\begin{align*}
    y &\geq 0, \\
    y &\geq \sum_{i = 1}^{n} w_i x_i + b, \\
    y &\leq \sum_{i = 1}^{n} w_i x_i + b - l(1-z), \\
    y &\geq uz.
\end{align*}

Big-M formulations are prone to numerical issues and increased running time when the bounds on the input are large, and thus benefit greatly from tighter bounds \cite{camm1990cutting,klotz2013practical}. For this reason, when optimising MIPs with embedded neural networks, either feasibility based or optimality based bound tightening is typically performed prior to optimisation \cite{grimstad2019relu,ceccon2022omlt}. It is therefore important that the user specifies bounds on the variables of the wider optimisation problem. Since this is not always possible, \scipml{} gives the user the option to decide between the classical big-M formulation or a novel SOS1-based formulation to represent ReLU activation functions. By doing the latter, we aim to maintain a valid formulation even in the absence of variable bounds, and minimise the pegging of solver performance to bound tightening procedures in the case of provided bounds. The downside of the SOS1-based formulation is that the SOS-constraint does not feature in the LP relaxation, and must be enforced by branching \cite{fischer2018branch}.

Our SOS1 formulation for an intermediate neuron with ReLU activation function introduces $s \in \mathbb{R}_{\geq 0}$ as a slack variable for the neuron.  The formulation is defined as:
\begin{align*}
    &y = \sum_{i = 1}^{n} w_i x_i + b + s, \\
    &\text{SOS1}(s, y).
\end{align*}

An SOS1 constraint enforces that at most one variable in the constraint can be non-zero. The above formulation ensures that the slack is only non-zero when $\sum_{i = 1}^{n} w_i x_i + b$ is negative. In this case, the slack must also take exactly the value of $-b - \sum_{i = 1}^{n} w_i x_i$, as the SOS1 constraint then forces $y$ to 0.

To illustrate the advantages of the SOS1 formulation, we present an example comparing it to the standard big-M formulation. Assume we have $l=-10^9$, $u=10^9$, and that the expression $\sum_{i = 1}^{n} w_i x_i + b$ evaluates to 0. We can consider the case that $z=10^{-6}$, which is within its integer feasibility tolerance of 0, and would be valid from the solver's perspective. Consider the constraints however:
\begin{align*}
    y &\geq 0 \\
    y &\geq 0 \\
    y &\leq 0 + 10^9 (1-10^{-6}) \\
    y &\leq 10^9 * 10^{-6} 
\end{align*}

In the big-M formulation for this poorly scaled example, the output $y$ is in the range $[-\epsilon, 10^3+\epsilon]$, where $\epsilon \in [0, 10^{-6}]$ is the constraint feasibility tolerance of the MIP solver. Thus, the output of the neural network for the MIP can vary greatly compared to its actual prediction from the same input. Meanwhile, due to our SOS1 formulation containing no big-M constraints, the error of the individual node's prediction is bounded by the feasibility tolerance value. We note however, that in practice such an occurrence is uncommon, since it occurs only when the LP solution happens to be integer feasible and the binary variable is free. We observed that when solving MIPs with embedded ML predictors using SCIP \cite{scip}, most binary variables are resolved through branching, with the MIP solver fixing the binary variable and not allowing any epsilon changes to its value. 

While we have focused on the ReLU activation function formulation, we note that \scipml also supports other activation functions such as sigmoid, and tanh. As \scipml uses SCIP \cite{scip} as a solver, these formulations can be added directly as non-linear expressions, requiring no reformulation or linearisation.

\subsection{Tree-Based Predictors}

To model tree-based predictors like decision trees, gradient boosted trees, or random forests, \scipml{} uses the formulation first introduced in Gurobi Machine Learning \cite{gurobiml}. As opposed to quadratic formulations \cite{mistry2021mixed} and linear big-M formulations \cite{ammari2023}, the Gurobi Machine Learning formulation uses a series of indicator constraints for each leaf node ensuring that any input maps to the appropriate leaf node of the tree. In what follows, we will give an overview of how decision trees are modelled in \scipml{}. Note that the MIP formulations of gradient boosted trees and random forests are simply a linear combination of the individual decision tree formulations.

For a decision tree, we denote the set of leaves as $\mathcal{L} \subseteq [m]$ and the input to the tree as $x \in \mathbb{R}^{n}$, where $m \in \mathbb{N}$ is the number of nodes of the tree and $n \in \mathbb{N}$ is the number of features. Each leaf $l \in \mathcal{L}$ is defined by a number of constraints on the input features of the tree that correspond to the branches taken in the path leading to $l$. We formulate decision trees by introducing one binary decision variable $\delta_l \in \{0,1\}$ for each leaf of the tree, which indicates whether the given input maps to $l$ (i.e., $\delta_l = 1$) or not (i.e., $\delta_l = 0$). Since in a decision tree exactly one leaf is chosen, we have to add the following constraint: 
\[
    \sum_{l \in \mathcal{L}} \delta_l = 1.
\]
To ensure that the input vector maps to the correct leaf, however, we need to introduce additional notation and constraints. For a node $v \in [m] \setminus \mathcal{L}$ in the tree, we denote by $i_v \in [n]$ the feature used for splitting, i.e., making the decision, and by $\theta_v \in \mathbb{R}$ the value at which the split is made. At a leaf $l \in \mathcal{L}$, we have a set $\mathcal{L}_l$ of inequalities of the form $x_{i_v} \leq \theta_v$ corresponding to the left branches leading to $l$ and a set $\mathcal{R}_l$ of inequalities of the form $x_{i_v} > \theta_v$ corresponding to the right branches. For each leaf, the inequalities describing $\mathcal{L}_l$ and $\mathcal{R}_l$ are imposed with indicator constraints modelling the following relationships: 
\begin{align*}
    \delta_l = 1 \implies &x_{i_v} \leq \theta_v - \frac{\epsilon}{2}, \quad \forall x_{i_v} \leq \theta_v \in \mathcal{L}_l, \\
    \delta_l = 1 \implies &x_{i_v} \geq \theta_v + \frac{\epsilon}{2}, \quad \forall x_{i_v} > \theta_v \in \mathcal{R}_l.
\end{align*}
In our implementation, $\epsilon \geq 0$ can be specified by a keyword parameter \texttt{epsilon} in the central function $\texttt{add\_predictor\_constr}$ for embedding a trained predictor as described in Section \ref{section:package}. By default, \scipml{} sets $\epsilon = 0$. Note that when $\epsilon$ is smaller than the default feasibility tolerance in SCIP\footnote{In SCIP 9.0, the feasibility tolerance is set to $10^{-6}$.} (as it is by default), and you have a solution where $| x_{i_v} - \theta_v | \approx \epsilon$, then SCIP can select an arbitrary child node of that decision in the tree. Note also that larger values for $\epsilon$ risk removing solutions containing values $x_{i_v} \approx \theta_v$. These removed solutions can change the optimal solution of the instance, or make the instance infeasible. 

\subsection{Argmax Formulation}

For many classification tasks, the ML predictor's evaluation is determined by an argmax function call. That is, if each class $j$ from the index set $\mathcal{J}$ has some predicted score $y_j \in \mathbb{R}$, the predictor outputs the class with the highest score. To model this functionality with MIP, we define a variable $m \in \mathbb{R}$, binary variables $z_j \in \{0,1\}, j \in \mathcal{J}$, and slack variables $s_j \geq 0, j \in \mathcal{J}$. The formulation we use is then given by:
\begin{align*}
    &y_j + s_j - m = 0, \ \forall j \in \mathcal{J} \\
    &\text{SOS1}(z_j, s_j), \ \forall j \in \mathcal{J} \\
    &\sum_{j = 1}^{|\mathcal{J}|} z_j = 1. 
\end{align*}

In the above formulation, the variable $m$ implicitly represents max$\{y_j: j \in \mathcal{J}\}$: If $m$ is less than the maximum, there exists a $j' \in \mathcal{J}$ such that $y_{j'} \geq m$, and for which $y_{j'} + s_{j'} - m = 0$ cannot be valid. On the other hand, if $m$ is greater than the maximum, then all slack variables are non-zero, and $\sum_{i = 1}^{|\mathcal{J}|} z_j = 1$ cannot be satisfied due to the SOS1 constraints. Because $m$ is the maximum, exactly one $z_{j'}$ is set to 1, where $y_{j'} = \text{max}\{y_j: j \in \mathcal{J}\}$. 

To truly return the argmax instead of a binary variable corresponding to the location of the maximum value, one can introduce a variable $a \in \mathbb{Z}$, and the corresponding constraint
\begin{align*}
    a = \sum_{j = 1}^{|\mathcal{J}|} j * z_j.
\end{align*}
% The variable $a$ will then have value corresponding to that returned by the traditional argmax function.

\section{Instance Library}
\label{section:library}

Accompanying \scipml{}, we present a new instance library, SurrogateLIB
\footnote{\url{https://zenodo.org/records/10357875}},
%\footnote{A compressed version of the library is provided in the supplementary materials. To respect the size limit of the submission, we had to remove instances.},
consisting of MIPs that contain embedded ML predictors as either constraints or components of the objective function.

SurrogateLIB is motivated by the large availability of data in the ML community and the need for more homogeneous non-trivial instance sets in the MIP community. Its goal is to be a library of instances with controllable complexity that is representative of MIPs with embedded ML predictors. To achieve this, semi-realistic MIP scenario generators are constructed in an array of application areas, including auto manufacturing, water treatment, and adversarial attacks, where real-world data is readily available for training the embedded ML predictors. By avoiding synthetically generated data to train the ML predictors, we aim to showcase the practical benefits of embedding predictive models in optimisation problems. An overview of each problem currently contained in SurrogateLIB, including detailed descriptions of the problems and how they are modelled, can be found in Appendix~\ref{appendix}. The current iteration of SurrogateLIB contains 1016 instances.

Alongside the library, we also include instance generators of the introduced problems that allow users to create homogeneous instances of varying difficulty. All generators feature two random seeds as common arguments. One is for randomising the ML predictor's training, and the other is for randomising some constraints of the MIP formulation. Thus, \scipml{} and SurrogateLIB present an opportunity to generate semi-realistic MIPs of controllable complexity by adjusting the size of the embedded ML predictors as well as the parameters of the instances. This way, we provide a set of homogeneous instances that are sufficiently diverse and relatively easy to solve without being trivial. Table \ref{tab:surrogateLIB} presents a summary for the current generators in SurrogateLIB. 

\begin{table}[t]
\centering
\resizebox{\columnwidth}{!}{%
\begin{tabular}{l|cccc}
\hline
\hline
\\[0.05cm]
Example Origin & Output Type & \#(Fixed) Input & \# ML Pred. & Brief Description \\[0.3cm]
\hline
\hline
\\
%Adversarial Attack \cite{lecun1998gradient,goodfellow2014explaining}
\begin{tabular}{@{}l@{}} Adversarial\\Attack \cite{lecun1998gradient,goodfellow2014explaining}\end{tabular}
& MR & (0)/256 & 1 & 
%Maximise chance of incorrect image classification \\
\begin{tabular}{@{}c@{}} Maximise chance of incorrect\\image classification\end{tabular} \\
&&&& \\

Auto Manufacturer \tablefootnote{\url{https://www.kaggle.com/datasets/gagandeep16/car-sales}} & R & (0)/10 & 3 & 
%Maximise sales of sufficiently different car \\
\begin{tabular}{@{}c@{}} Maximise sales of\\sufficiently different car\end{tabular} \\
&&&& \\

City Planning \tablefootnote{\url{https://www.kaggle.com/datasets/krzysztofjamroz/apartment-prices-in-poland}} & R & (7)/14 & 25 &
%Maximise living quality of new residents \\
\begin{tabular}{@{}c@{}} Maximise living quality\\of new residents\end{tabular} \\
&&&& \\

%Function Approximation \cite{grimstad2019relu}
\begin{tabular}{@{}l@{}} Function\\Approximation \cite{grimstad2019relu}\end{tabular}
& R & (0)/5 & 2 &
%Maximise function while satisfying equality constraints \\
\begin{tabular}{@{}c@{}} Maximise function while\\satisfying equality constraints\end{tabular} \\
&&&& \\

%Patalable Diet Problem \cite{peters2021nutritious,opticl}
\begin{tabular}{@{}l@{}} Patalable Diet\\Problem \cite{peters2021nutritious,opticl}\end{tabular}
& R & (2)/25 & 1 & 
\begin{tabular}{@{}c@{}} Minimise transport cost while meeting\\nutritional / palatability requirements\end{tabular} \\
&&&& \\
%Minimise transport cost while meeting nutritional / palatability requirements \\

Tree Planting \cite{wood2023tree} & R & (6)/7 & 100 & 
Maximise tree survival rates \\
&&&& \\

Water Potability \tablefootnote{\url{https://github.com/MainakRepositor/Datasets/tree/master}} & C & (0)/9 & 50 &
%Maximise amount of drinkable water after treatment \\
\begin{tabular}{@{}c@{}} Maximise amount of drinkable\\water after treatment\end{tabular} \\
&&&& \\

Wine Quality \cite{cortez2009modeling} & R & (0)/11 & 5 & 
%Maximise quality of mixed wine from suppliers \\
\begin{tabular}{@{}c@{}} Maximise quality of mixed\\wine from suppliers\end{tabular} \\
&&&& \\

%Workload Dispatching\cite{eml}
\begin{tabular}{@{}l@{}} Workload\\Dispatching\cite{eml}\end{tabular}
& R & (0)/3 & 48 &
%Maximise minimum efficiency of compute cores \\
\begin{tabular}{@{}c@{}} Maximise minimum efficiency\\of compute cores\end{tabular} \\[0.5cm]
\hline
\hline
\end{tabular}
}
\caption{Overview of instance generators available in SurrogateLIB. \emph{Output type} is labelled as (M)ulti-(R)regression or (C)lassification. \emph{\#(Fixed) Input} refers to the number of MIP variables that are input to each embedded ML predictor, both for those that can take free values and for those that are fixed. \emph{\#ML Pred.} is the number of embedded ML predictors by default in the generator.}
\label{tab:surrogateLIB}
\end{table}

\section{Computational Experiments}
\label{section:experiments}

One question that often arises when practitioners deal with surrogate models is the following: How large of a ML predictor can actually be embedded without blowing up solution time? To answer this question, we will present a set of experiments aiming to provide some intuition for the scale of embedded ML predictors that practitioners can expect to optimise over.

To determine the impact of embedding a ML predictor on the solution time, we use SurrogateLIB and progressively embed larger ML predictors until the solution time is 30 minutes or greater. We embed fully connected feed forward neural networks with ReLU activation functions and gradient-boosted decision trees (GBDT). For the neural network, we test both the SOS1 and big-M formulation, where the hidden layer size is fixed to 16. For GBDTs we fix the maximum depth of each estimator to 5. Increasing the size of each ML predictor then corresponds to either introducing an additional hidden layer or an additional estimator. The function approximation problem has been removed so as to leave only semi-real-world examples. GBDT is not used for problems where the embedded ML predictors have fewer than $4$ non-fixed inputs, and for the image based adversarial attack problem. 

We fixed the hidden layer sizes of the neural networks and depth of individual estimators in the GDBTs as we wanted only a single axis on which to increase the size of the embedded predictors. Size 16 was chosen as it is a common choice for smaller neural network architectures, and depth 5 was chosen as it is appropriate for a wide range of SurrogateLIB instances. We reduced the experiment to only two types of ML predictors as they admit valid integer linear formulations, can be scaled, and encompass simpler ML predictors, e.g., decision trees. 

For all experiments, SCIP 9.0.0 \cite{scip9} is used for MIP solving with SOPLEX 7.0.0 \cite{scip9} as the LP solver, PyTorch 2.1.0 \cite{paszke2019} is used for training the neural networks, and Scikit-learn 1.4.2 \cite{scikit-learn} is used for training the GBDTs. All experiments are run on a cluster equipped with Intel Xeon Gold 5122 CPUs with 3.60GHz and 96GB main memory. 

During preliminary experiments, we observed a common problem where deeper neural networks would collapse into predicting a constant output. We believe this is caused by vanishing gradients, and is the result of using too deep of a neural network to represent simpler relations. The resulting MIP with the collapsed neural network embedded is often quickly identified as infeasible, as the constant output cannot be changed to satisfy some constraints. Such a result is ultimately misleading for determining the scale of ML predictors that can be embedded. To overcome this issue, we initialised the neural networks with a Xavier uniform distribution \cite{glorot2010understanding}, introduced gradient clipping \cite{goodfellow2016deep}, and added an L2 penalty. Additionally, we stopped the experiment at a maximum of 15 layers or estimators. 

\begin{table}[t]
\centering
\footnotesize
\renewcommand{\arraystretch}{1.7} % Adjust the space between rows
\begin{tabular}{l|c|c|c|c}
\hline
\hline
 & Adversarial Attack & Auto Manufacturer & City Planning & Platatable Diet \\
\hline
\hline
ReLU SOS & 6 & 2 & 0 & 14* \\
%\hline
ReLU Big-M & 6 & 5 & 0 & 13* \\
%\hline
GBDT & - & 15+ & 2 & 15+ \\
%\hline
\hline
\hline
 & Tree Planting & Water Potability & Wine Manufacturer & Workload Dispatching \\
\hline
\hline
ReLU SOS & 0 & 8 & 3 & 0 \\
%\hline
ReLU Big-M & 2* & 8 & 3 & 0 \\
%\hline
GBDT & - & 6 & 4 & - \\
\hline
\hline
\end{tabular}
\caption{Maximum embedding size for ML predictors in SurrogateLIB such that the resulting MIP can be solved within 30 minutes. Entries with "0" indicate that the starting embedding size could not be solved for the corresponding problem-formulation pair. Entries with "-" indicate that the corresponding experiment was skipped. \\
\footnotesize{*Unreliable results due to collapsing neural network.}
}
\label{tab:experiments-surrogatelib}
\end{table}

The summarised result for the experiment is presented in Table \ref{tab:experiments-surrogatelib}. The table shows that the scale of ML predictors that can be practically embedded is highly predictor-dependent. For example, in the auto manufacturer instance, a predictor of scale 2, 5, or 15+ can be solved to optimality in under 30 minutes depending on the embedded predictor type and formulation. Furthermore, we can also observe a dependence on the instance-type. Whereas \scipml{} is able to solve instances of the Water Potability problem with relatively large embedded ML predictors, it fails to solve the Workload Dispatching problem with the smallest embedded ML predictor size. This is despite the Water Potability problem having a larger input space for each embedded ML predictor, and more embedded ML predictors. This shows that SurrogateLIB contains a broad range of generators of varying difficulty, making it a versatile tool to test new approaches for solving MIPs with surrogate models.

\section{Conclusion}

In this paper we have introduced the Python package \scipml, which automatically embeds trained ML predictors in MIPs. We have briefly introduced the API, expanded on some of the core MIP formulations including novel formulations for ReLU neural networks and the argmax function, and presented a set of computation results that provide an intuition for the scale of ML predictors that can be embedded in practice. For the final point, we introduced a library of homogeneous semi-real-world MIP instance generators, which form the extendable library SurrogateLIB. 

\section*{Acknowledgements}
We thank Mohammed Ghannam, Robert Luce, Julian Manns, and Franziska Schlösser for their help with deploying SCIP via PyPI. The work for this article has been conducted in the Research Campus MODAL funded by the German Federal Ministry of Education and Research (BMBF) (fund numbers 05M14ZAM, 05M20ZBM).

\ifarxiv
\bibliographystyle{unsrt}
\else
\bibliographystyle{splncs04}
\fi
\bibliography{mybib}

\newpage

\appendix

\section{SurrogateLIB}
\label{appendix}

A common problem when formulating a MIP formulation of a system is the presence of unknown or highly complex relationships. To address this challenge, ML predictors are often used to approximate these relationships and are then embedded in the MIP as surrogate models. Derived from this terminology, we present \emph{SurrogateLIB}, a new library of MIP instances. 

SurrogateLIB is motivated by the large availability of data in the ML community and the need for more homogeneous non-trivial instance sets in the MIP community. Its goal is to be a library of instances with controllable complexity that is representative of MIPs with embedded ML predictors. To achieve this, semi-realistic MIP scenario generators are constructed in an array of application areas where real-world data is available for training the embedded ML predictors. Each generator features two random seeds as common arguments, which control the ML predictor initialisation and some randomisation of the MIP budget constraints.

The naming convention of an instance from SurrogateLIB requires us to introduce some notation. Let $A$ be the problem type, $B$ be the problem parameters, $C$ be the ML predictor type, $D$ be the predictor parameters, $E$ be the ML framework, and $F$, $G$ be the data and training random seeds. An instance from SurrogateLIB is then named:

\begin{center}
    $A\_B\_C\_D\_E\_F\_G.mps$
\end{center}

We explicitly list the problem types ($A$) and the associated problem parameters ($B$) in Table \ref{tab:problem_type_surrogatelib}, where thorough definitions of each generator are given in the subsections of the Appendix. 

\begin{table}[t]
\centering
\footnotesize
\renewcommand{\arraystretch}{1.3} % Adjust the space between rows
\begin{tabular}{l|c|c|c|c|c|c|c|c|c}
$A$ & adversarial & auto & city & function & palatable & tree & water & wine & workload \\
\hline
$B$ & n & - & n & n & - & n & n & n, m & - 
\end{tabular}
\caption{Simplified naming scheme of each problem type and the problem parameters as detailed in the full descriptions}
\label{tab:problem_type_surrogatelib}
\end{table}

The predictor types ($C$) that we have currently embedded as part of SurrogateLIB are: linear (linear regression), dt (decision tree), gbdt (gradient boosted decision tree), rf (random forest), mlp-sos / mlp-bigm (feed forward neural network with ReLU activation functions and the SOS1 / Big-M formulation respectively). The parameters ($D$) for each predictor that follow in the naming scheme are given in Table \ref{tab:predictor_type_surrogatelib}.

\begin{table}[t]
\centering
\footnotesize
\renewcommand{\arraystretch}{1.3} % Adjust the space between rows
\resizebox{\columnwidth}{!}{%
\begin{tabular}{l|c|c|c|c|c|c|c|c|c}
$C$ & linear & dt & gbdt & rf & mlp-sos & mlp-bigm\\
\hline
$D$ & - & depth & \# estimators, depth & \# estimators, depth & \# layers, layer size & \# layers, layer size
\end{tabular}
}
\caption{ML predictor parameters for SurrogateLIB}
\label{tab:predictor_type_surrogatelib}
\end{table}

The list of ML frameworks, i.e., valid values of $E$, consists of: sk (scikit-learn \cite{scikit-learn}), torch (PyTorch \cite{paszke2019}), keras (Keras \cite{keras}), lgb (LightGBM \cite{ke2017lightgbm}, and xgb (XGBoost \cite{xgboost}). The data and training seed, $F$ and $G$, can independently take the values $0$ and $1$.

We now present each generator for SurrogateLIB in the following subsections.

\subsection{Adversarial Example}

An adversarial example is a purposefully altered input instance to a trained ML predictor such that its resulting label is incorrect. Adversarial examples are of extreme importance to many applications \cite{goodfellow2014explaining}, with an entire field devoted to robust ML training \cite{wong2018provable}. They feature as a prominent example for embedded ML predictors in MIPs as they require only the ML predictor and an example input to construct the entire MIP. We use the well-known MNIST dataset \cite{lecun1998gradient} in our generator of adversarial example MIPs. 

Let $f:[0,1]^{n\times n} \rightarrow \reals^{10}$ be the trained ML predictor that takes as input an $n\times n$ image of a black and white photo and outputs the probability of each digit $[0,9]$ being the image, where $n \in \mathbb{N}$. The variables $\bar{x}_{i,j} \in [0,1]$ for $i,j \in [n]$ are the pixel values of the input instance after perturbation, where $x_{i,j} \in [0,1]$ for $i,j \in [n]$ is the given image. The variables $d_{i,j} \in [0, 1]$ for $i,j \in [n]$ represent the difference between $\bar{x}_{i,j}$ and $x_{i,j}$, and $\delta \in \reals_{\geq 0}$ represents the total difference budget. Finally, $y \in \reals^{10}$ are the output probabilities, where without loss of generality $y_{1}$ is the true label and $y_{2}$ is the largest valued incorrect label for the base image. The MIP is defined as:
\begin{equation*}
\begin{aligned}
& \underset{x}{\text{min}} & &y_2 - y_1 \\
& \text{s.t.} & &\bar{x}_{i,j} - x_{i,j} \leq d_{i,j}, \;\;\forall i,j \in [n] \\
&&& x_{i,j} - \bar{x}_{i,j} \leq d_{i,j}, \;\;\forall i,j \in [n] \\
&&& \sum_{i=1}^{n} \sum_{j=1}^{n} d_{i,j} \leq \delta \\
&&& y = f(x).
\end{aligned}
\end{equation*}

\subsection{Auto Manufacturer}

For this generator, we take the point of view of an auto manufacturer that needs to design a new vehicle. Based on historical data about past sales, the manufacturer aims to maximise the number of vehicles sold, while ensuring that the designed vehicle is sufficiently different to any other popular vehicle already on the market. Additionally, the manufacturer wants to ensure that the designed vehicle has a high resell value and is fuel-efficient.

Let $x \in \mathbb{R}^n$ be the vector of variable vehicle features that the manufacturer can control, where for our case $n = 10$. These features are: vehicle type, engine size, horsepower, wheelbase, width, length, curb weight, fuel capacity, efficiency and power performance factor. Every feature has a lower bound $l_i \in \reals$ and an upper bound $u_i \in \reals$, that is, $l \leq x \leq u$ with $l, u \in \mathbb{R}^n$. Furthermore, let $y_{\text{price}} \in \reals$, $y_{\text{resell}} \in \reals$, and $y_{\text{sold}} \in \reals$ be the features of the vehicle that the manufacturer cannot control. In this example, these uncontrollable features are the price of the vehicle, its resell value, and the number of sold vehicles.

Since the relationship between $x$ and $(y_{\text{price}},y_{\text{resell}},y_{\text{sold}})$ is unknown, the auto manufacturer uses ML to learn from real-world historic data\footnote{\url{https://www.kaggle.com/datasets/gagandeep16/car-sales}}. To do so, a ML predictor, $f: \reals^{10} \rightarrow \reals$, is trained that takes as input the controllable features of the vehicle and outputs the predicted price. Two more ML predictors, $g: \reals^{11} \rightarrow \reals$ and $h: \reals^{11} \rightarrow \reals$, are trained, which take as input the controllable features and the predicted price of the vehicle, and output the predicted amount of vehicles sold and predicted resell value respectively.

Recall that the manufacturer aims to design a vehicle different to what already exists, requiring us to introduce additional constraints. For a finite index set $\mathcal{V}$, we denote by $\bar{x}^j \in \mathbb{R}^{10}, j \in \mathcal{V}$, the controllable features of the $|\mathcal{V}|$ many most popular vehicles. To ensure that $x$ is sufficiently different to each $\bar{x}^j$, we consider the relative difference of each feature and enforce that the cumulative difference among all features are above a given threshold $\gamma \in \reals$. Formally, we obtain the constraints
\begin{align}
    \sum_{i = 1}^{10} \frac{| x_i - \bar{x}_i^j |}{u_i - l_i} \geq \gamma, \quad \forall j \in \mathcal{V}. \label{eq:diff_cons}
\end{align}
Note that even though the absolute value function is non-linear, it can be modelled using linear constraints and binary variables. Finally, let $\alpha \in [0,1]$ be the minimum ratio of the resell value compared to the actual purchase price, and let $\beta \in \reals$ be the minimum level of fuel-efficiency, where without loss of generality $x_0$ is the feature index for fuel-efficiency. The MIP is defined as:

\begin{equation*}
\begin{aligned}
& \underset{x}{\text{max}} & & y_{\text{sold}} \\
& \text{s.t.} & & y_{\text{price}} = f(x) \\
&&& y_{\text{sold}} = g(x, y_{\text{price}}) \\
&&& y_{\text{resell}} = h(x, y_{\text{price}}) \\
&&& y_{\text{resell}} \geq \alpha y_{\text{price}} \\
&&& x_{0} \geq \beta \\
&&& \eqref{eq:diff_cons} \\
&&& l \leq x \leq u.
\end{aligned}
\end{equation*}

\subsection{City Manager}

For this generator, we take the point of view of a city council that has been tasked with building amenities missing in the community. These amenities include two schools, two doctor clinics, two post-offices, two kindergartens, two restaurants, two pharmacies, and a college. As we do not have a direct measure for quality of life improvements by the distance to these amenities, we use the price of apartments as a proxy, for which data is readily available \footnote{\url{https://www.kaggle.com/datasets/krzysztofjamroz/apartment-prices-in-poland}}. 

Let $f:\reals^{14} \rightarrow \reals$ be the trained ML predictor that takes as input a feature vector of an apartment and outputs the predicted price. A feature vector $x \in \reals^{14}$ of an apartment consists of the square meters, number of rooms, floor level, whether the apartment has parking, a balcony, or security, and the shortest distance to the town center, a doctor's office, a post-office, a kindergarten, a restaurant, a pharmacy, and a college. We construct an imaginary town of $n \in \mathbb{N}$ many apartments in an $m \times m$-km grid, where $m \in \reals_{\geq 0}$. The features of each apartment that are not distance based are randomised and fixed. This includes the location of the apartment, and therefore the distance to the town center. The output of the ML predictor for apartment $i \in [n]$ is represented by the variable $y_{i} \in \reals$.

Let $\mathcal{J} = \{\text{doctor, post, kinder, restaurant, pharmacy, college}\}$ be the set of amenities. For each amenity $j \in \mathcal{J}$, let $b^{j}_{k,0} \in \reals$ be the variables that represent the grid-coordinate of the first amenity of that type constructed for dimension $k \in \{1,2\}$ of the grid. The variables $b^{j}_{k,1} \in \reals$ then represent the $k \in \{1,2\}$ dimension grid-coordinate of the second amenity of that type for all amenities except for college. For the randomly generated apartments, let $a_{i,k} \in \reals$ be the $k \in \{1,2\}$ dimension grid-coordinate of apartment $i \in [n]$. The feature vector $x^{i} \in \reals^{14}$ then represents apartment $i \in [n]$. Additionally, let $\gamma \in \reals_{\geq 0}$ be a minimum distance requirement between any two amenities of the same type. The MIP is defined as:

\begin{equation*}
\begin{aligned}
& \underset{x}{\text{max}} & & \sum_{i=1}^{n} y_{i} \\
& \text{s.t.} & & y_{i} = f(x^{i}), \;\; \forall i \in [n] \\
&&& \sum_{k =1}^{2} |b^{j}_{k,0} - b^{j}_{k,1} | \geq \gamma, \;\; \forall j \in \mathcal{J} \setminus \{\text{college}\} \\
&&& x^{i}_{j} = \text{min}( \sum_{k =1}^{2} |b^{j}_{k,0} - a_{i,k} |, \sum_{k =1}^{2} |b^{j}_{k,1} - a_{i,k} |) \;\; \forall j \in \mathcal{J} \setminus \{\text{college}\}, i \in [n] \\
&&& x^{i}_{\text{college}} = \sum_{k =1}^{2} |b^{\text{college}}_{k,0} - a_{i,k} | \;\; \forall i \in [n] 
\end{aligned}
\end{equation*}

Note that even though the absolute value and min functions are non-linear, they can be modelled via linear constraints and additional binary variables.

\subsection{Function Approximation}

As opposed to all other members of SurrogateLIB, this MIP is entirely contrived from random data. It was the MIP used in \cite{grimstad2019relu} to evaluate the impact of different big-M formulations and bound-tightening procedures on solver performance. 

Let $f: \reals^{n} \rightarrow \reals$ and $g: \reals^{n} \rightarrow \reals$ be trained ML predictors, where $n \in \mathbb{N}$. Each predictor approximates a different random quadratic function of the form $x^{\intercal}Qx + a^{\intercal}x + b$, where $x \in \reals^{n}$, $Q \in \reals^{n \times n}$, $a \in \reals^{n}$, and $b \in \reals$. The values of $Q, a$, and $b$ are samples from a uniform distribution. Additionally, a constant value $c \in \reals$ is generated. The MIP is defined as:

\begin{align*}
& \underset{x}{\text{min}} \;\; f(x) \\
& \text{s.t.} \;\;\; g(x) = c
\end{align*}

\subsection{Palatable Diet}

For this generator, we take the view point of the World Food Programme that must decide the quantities of different food types that it will deliver. This scenario phrased as an optimisation problem has been discussed heavily in \cite{peters2021nutritious}, and was a leading example for embedded ML predictors in MIP for the framework OptiCL \cite{opticl}. Borrowing data that they have made publicly available, we construct a simplified variant of the original food delivery problem, which was coined \textit{The Palatable Diet Problem}. The name is derived by extending the traditional diet problem via a minimum palatability constraint.

Let $x \in \reals^{25}_{\geq 0}$ be a feature vector that represents a basket of 25 types of food. Additionally, let $n_{i,j} \in \reals_{\geq 0}$ be the nutritional value of nutrient $j \in \mathcal{J}$ for one unit of food $i \in [25]$, where for our scenario we consider 12 nutrients, i.e., $|\mathcal{J}|=12$. The value $\gamma_{j} \in \reals_{\geq 0}$ is a minimum requirement of nutrient $j \in \mathcal{J}$. While some baskets of food would meet nutritional requirements, they are undesirable from a consumers perspective. Therefore, a palatability requirement is enforced, with this requirement being derived from data due to the complexity of defining such a constraint. Let $f: \reals^{25} \rightarrow \reals$ be a trained ML predictor that predicts the palatability of a given food basket, and $\beta \in \reals$ be some minimum requirement. Additionally, let $c_{i} \in \reals_{\geq 0}$ be the transportation cost of a single unit of food for $i \in [25]$. The MIP is defined as:

\begin{align*}
& \underset{x}{\text{min}} \;\; \sum_{i=1}^{25} c_{i} x_{i} \\
& \text{s.t.} \;\;\; \sum_{i=1}^{25} n_{i,j} x_{i} \geq \gamma_{j}, \; \forall j \in \mathcal{J} \\
& \;\;\;\;\;\;\; f(x) \geq \beta \\
& \;\;\;\;\;\;\; x_{\text{salt}} = 5 \\
& \;\;\;\;\;\;\; x_{\text{sugar}} = 20 
\end{align*}

\subsection{Tree Planting}

For this generator, we take the point of view of an individual responsible for replanting a set of trees on a stretch of empty land. The goal of this problem is to plant a diverse set of trees, while respecting a budget constraint, and maximising the total amount of surviving trees. Thanks to studies on the effect of soil sterilisation and light levels on tree survival rates \cite{wood2023tree}, we have data readily available for such a problem. 

Let $x_{i} \in \reals^{7}$ be a feature vector of planting location $i \in [n]$, where $n \in \mathbb{N}$. We randomly generate a grid of planting locations with randomised feature values, where the binary decision of whether the location is sterilised set to false. A sterilisaion budget is then introduced and denoted by $\beta \in \mathbb{N}$. For each planting location we must plant one of four tree species $k \in [4]$. Each species has a minimum expected survival requirement $\gamma_{k} \in \reals_{\geq 0}$, and a planting cost $c_{k} \in \reals_{\geq 0}$. For each species, a trained ML predictor $f_{k}: \reals^{7} \rightarrow \reals$ is constructed that predicts the probability that a tree of species $k \in [4]$ survives in a given planting location. To model the MIP we introduce variables $s_{i,k} \in \reals$ representing the (adjusted) predicted survival rate of tree species $k$ in location $i$, and $s'_{i,k}$ to represent the adjusted survival rates. Binary variables $p_{ik} \in \{0,1\}$ represent whether tree species $k \in [4]$ was planted in location $i \in [n]$. The MIP is defined as:

\begin{align*}
& \underset{x}{\text{max}} \;\; \sum_{i=1}^{n} \sum_{k=1}^{4} s'_{i,k} \\
& \text{s.t.} \;\;\; \sum_{i=1}^{n} s'_{i,k} \geq \gamma_{j}, \; \forall k \in [4] \\
& \;\;\;\;\;\;\; p_{i,k} = 1 \xrightarrow{} s'_{i,k} \leq s_{i,k}, \; \forall i \in [n], k \in [4] \\
& \;\;\;\;\;\;\; p_{i,k} = 0 \xrightarrow{} s'_{i,k} \leq 0, \; \forall i \in [n], k \in [4] \\
& \;\;\;\;\;\;\; \sum_{k=1}^{4} c_{k} p_{i,k}, \; \forall i \in [n] \\
& \;\;\;\;\;\;\; \sum_{k=1}^{4} p_{i,k} = 1, \; \forall i \in [n] \\
& \;\;\;\;\;\;\; \sum_{i=1}^{n} x^{\text{sterile}}_{i} \leq \beta, \\
& \;\;\;\;\;\;\; f_{k}(x_{i}) = s_{i,k}, \; \forall i \in [n], k \in [4]
\end{align*}

\subsection{Water Potability}

For this generator, we take the point of view of a health organisation that is treating drinking water. The organisation has access to processes that decrease or increase levels of attributes for a given sample of water. There are however limitations on the total amount that each attribute can be increased or decreased. The goal of the optimisation problem is to maximise the amount of drinkable water given a set of non-drinkable water. Data for this problem is publicly available \footnote{\url{https://github.com/MainakRepositor/Datasets/tree/master}}. 

Let $x_{i} \in \reals^{9}$ be the feature vector of given water sample, where we have taken $n \in \mathbb{N}$ many undrinkable water samples. For each water sample $i \in [n]$ and feature $j \in [9]$, let $a_{i,j} \in \reals_{\geq 0}$ be the positive change after applying treatment, and $b_{i,j} \in \reals_{\geq 0}$ be the negative change after applying treatment. For each feature $j \in [9]$ there is a positive and negative budget $\gamma^{+}_{j}, \gamma^{-}_{j} \in \reals_{\geq 0}$. Before treating the water, the undrinkable samples have values $w_{i} \in \reals^{9}$ for sample $i \in [n]$. To determine if the water is drinkable, we have access to a trained ML predictor $f: \reals^{9} \rightarrow \reals$. The ML predictor takes the feature vector of a water sample as input and outputs whether or not it is drinkable. The binary variables that determine if water sample $i \in [n]$ is drinkable are $y_{i} \in \{0,1\}$. The MIP is defined as:
\begin{equation*}
\begin{aligned}
& \underset{x}{\text{max}} & & \sum_{i=1}^{n} y_{i} \\
& \text{s.t.} & & x_{i,j} = w_{i,j} + a_{i,j} - b_{i,j}, \; \; \forall i \in [n], j \in [9] \\
&&& \sum_{i=1}^{n} a_{i,j} \leq \gamma^{+}_{j}, \; \; \forall j \in [9] \\
&&& \sum_{i=1}^{n} b_{i,j} \leq \gamma^{-}_{j}, \; \; \forall j \in [9] \\
&&& y_{i} = f(x_{i}), \; \; \forall i \in [n]
\end{aligned}
\end{equation*}

\subsection{Wine Blending}

For this generator, we take the point of view of a wine manufacturer. This manufacturer does not grow the grapes themselves, but rather purchases grapes from other growers and blends the result into a new wine. The grapes vendors sell grapes because a wine made from their grapes alone is predicted to be of poor quality. The objective of this problem is to create a bouquet of wines of fixed size with highest average quality. Each vineyard has an associated costs for purchasing a unit of grapes, and the wine manufacturer has a total budget. We assume that readily available data from features of wine \cite{cortez2009modeling} can be transferred to features of grapes for this optimisation problem, and that the features of a resulting wine are a linear combination of those in the blend.

Let $x_{i} \in \reals^{11}_{\geq 0}$ be the feature vector of a grape (wine) blend for $i \in [n]$, and $w_{j} \in \reals^{11}_{\geq 0}$ be the fixed feature vector of a wine (grape) from vendor $j \in [m]$. The features of a grape (wine) are fixed acidity, volatile acidity, citric acid, residual sugar, chlorides, free sulfur, total sulfur, density, pH, sulphates, and alcohol percentage. To determine wine quality, we have access to a trained ML predictor $f: \reals^{11} \rightarrow \reals$, which predicts the wine quality of the grape blend $i \in [n]$ from its associated feature vector $x_{i}$. The variables $y_{i} \in [0,10]$ for $i \in [n]$ represent the predicted quality, where $\alpha \in [0,1]$ is a minimum required quality.

Each grape vendor $j \in [m]$ has an amount of available grapes $\gamma_{j} \in \reals_{\geq 0}$, and an associated cost per unit $c_{j} \in \reals_{\geq 0}$. The variables $b_{i,j} \in \reals_{\geq 0}$ represent the amount of grapes used for wine blend $i \in [n]$ from vendor $j \in [m]$, and $\beta \in \reals_{\geq 0}$ is the total budget for the purchasing of all grapes. We assume that all created wine blends are of the same size. The MIP is defined as:

\begin{equation*}
\begin{aligned}
& \underset{x}{\text{max}} & & \sum_{i=1}^{n} y_{i} \\
& \text{s.t.} & & x^{i} = \sum_{j=1}^{m} b_{i,j} w_{j}, \; \; \forall i \in [n] \\
&&& \sum_{i=1}^{n} b_{i,j} \leq \gamma_{j}, \; \; \forall j \in [m] \\
&&& \sum_{i=1}^{n} \sum_{j=1}^{m} c_{j} b_{i,j} \leq \beta, \\
&&& \sum_{j=1}^{m} b_{i,j} = 1, \; \; \forall i \in [n] \\
&&& y_{i} \geq \alpha, \; \; \forall i \in [n] \\
&&& y_{i} = f(x_{i}), \; \; \forall i \in [n]
\end{aligned}
\end{equation*}

\subsection{Workload Dispatching}

For this generator, we take the point of view of a job scheduler that must issue compute jobs to individual compute cores. The objective of the optimisation problem is to maximise the minimum efficiency of all cores. Similar to the palatable diet problem for OptiCL \cite{opticl}, the \textit{Workload Dispatching Problem} served as the leading example for the EML framework \cite{eml}. Borrowing data that they have publicly available\footnote{\url{https://bitbucket.org/m_lombardi/eml-aij-2015-resources/src/master/training/}}, we construct instances of the problem. 

As efficiency of a compute core is difficult to determine apriori, a ML predictor is used. Let $x_{i} \in \reals^{3}$ be the feature vector of a set of jobs assigned to compute core $i \in [n]$ for $n \in \mathbb{N}$. The features are the average clock (equivalently cycles) per instruction (CPI) of jobs assigned to the compute core, the average CPI of neighbouring cores, and the average CPI of non-neighbouring cores. The neighbourhood of a core $i \in [n]$ is defined as $\mathcal{N}(i) \subseteq [n]$. To obtain the feature sets, jobs must first be distributed to each core, where an equal amount of jobs are assigned to each core. For this instance $m \in \mathbb{N}$ many jobs are randomly generated with varying CPIs $c_{j}$ for $j \in [m]$. The variables $x_{i,j} \in \{0,1\}$ represent whether job $j \in [m]$ is assigned to compute core $i \in [n]$. A trained ML predictor $f: \reals^{3} \rightarrow \reals$ predicts the efficiency $y_{i} \in \reals$ of compute core $i \in [n]$ after all cores are assigned jobs. Finally, the variable $e \in \reals$ represents the minimum efficiency over all cores. The MIP is defined as:

\begin{equation*}
\begin{aligned}
& \underset{x}{\text{max}} & & e \\
& \text{s.t.} & & \sum_{i=1}^{n} x_{i,j} = 1, \; \; \forall j \in [m] \\
&&& \sum_{j=1}^{m} x_{i,j} = \frac{m}{n}, \; \; \forall i \in [n] \\
&&& x_{i}^{\text{avg}} = \frac{n}{m} \sum_{j=1}^{m} c_{j} x_{i,j}, \; \; \forall i \in [n] \\
&&& x_{i}^{\text{avg-neigh}} = \frac{1}{|\mathcal{N}(i)|} \sum_{i' \in \mathcal{N}(i)} x_{i'}^{\text{avg}}, \; \; \forall i \in [n] \\
&&& x_{i}^{\text{avg-far}} = \frac{1}{n - |\mathcal{N}(i)|} \sum_{i' \in [n] \setminus \mathcal{N}(i)} x_{i'}^{\text{avg}}, \; \; \forall i \in [n] \\
&&& y_{i} = f(x_{i}^{\text{avg}}, x_{i}^{\text{avg-neigh}}, x_{i}^{\text{avg-far}}) \\
&&& e \leq y_{i}, \; \; \forall i \in [n]
\end{aligned}
\end{equation*}

%%%%%%%%%%%%%%%%%%%%%%%%%%%%%%%%%%%%%%%%%%%%%%%%%%%

\end{document}